\documentstyle[12pt]{article}
\begin{document}

\rightline{FTUV/96-4, IFIC/96-5}
\rightline{January 1996}
\rightline{q-alg/9604024}

\begin{center}
{\Large {\bf Differential calculus on `non-standard' 
($h$-deformed)  Minkowski  spaces}}\footnote{Lecture delivered at
the 1995 Quantum Groups and Quantum Spaces Banach conference, Warsaw. 
To appear in the proceedings, R. Budzynski and S.
Zakrzewski eds., Banach Centre publications. AMS numbers 20G45, 
20N99, 81R50, 53Z05.}
\end{center}

\begin{center}
{\large{J.A. de Azc\'{a}rraga$^{\dagger}$
and F. Rodenas$^{\ddagger}$}}
\end{center}

\noindent
{\small{\it $ \dagger$ Departamento de F\'{\i}sica Te\'{o}rica and IFIC,
Centro Mixto Universidad de Valencia-CSIC
E-46100-Burjassot (Valencia), Spain.}}

\noindent
{\small{\it $ \ddagger $ Departamento de Matem\'atica Aplicada,
 Universidad Polit\'ecnica de Valencia
E-46071 Valencia, Spain.}}

\vspace{1\baselineskip}

\begin{abstract}
The differential calculus on `non-standard'  $h$-Minkowski spaces is given.
In particular it is  shown  that,
for them, it is possible to introduce  coordinates and derivatives
which are simultaneously hermitian.
\end{abstract}

\section{Introduction}

\indent

We review first the properties of the two
deformed   Minkowski spaces \cite{JPA96}  associated with the 
`Jordanian' or `non-standard' $h$-deformation, $SL_h(2)$, of $SL(2,C)$
\cite{DEMI,EOW,KAR}. The $GL_h(2)$  deformation is 
defined as the associative algebra generated by the 
entries  $a,b,c,d$ of a matrix $M$, the commutation properties of which
may be  expressed by an `FRT'  \cite{FRT1} equation,
$R_{12}M_1M_2=M_2M_1R_{12}$, in which $R$  is
the triangular solution of the Yang-Baxter equation 
\begin{equation}\label{1.5}
R_h= \left[ 
\begin{array}{cccc}
1 & -h & h & h^2 \\
0 & 1 & 0 & -h \\
0 & 0 & 1 & h \\
0 & 0 & 0 & 1 
\end{array} \right] \;,
\; \hat{R}_h \equiv {\cal P}R_h= \left[ 
\begin{array}{cccc}
1 & -h & h & h^2 \\
0 & 0 & 1 & h \\
0 & 1 & 0 & -h \\
0 & 0 & 0 & 1 
\end{array} \right] \;,   \; {\cal P}R_h{\cal P}=R_h^{-1} \;.
\end{equation}	
\noindent
The commutation relations of the group algebra  generators in $M$ are 
\begin{equation}\label{1.6}
\begin{array}{lll}
{} [a,b]= h(\xi-a^2) \quad, & \quad [a,c]=hc^2 \quad,&\quad 
[a,d]= hc(d-a) \quad, \\
{}   [b,c]=h(ac+cd) \quad,  & \quad
 [b,d]=h(d^2-\xi) \quad, & \quad [c,d]=-hc^2 \quad ;
\end{array}
\end{equation}	
\begin{equation}\label{1.7}
\xi \equiv det_hM = ad - cb - h cd  \quad.
\end{equation}	
\noindent
Setting $\xi =1$  reduces  $GL_h(2)$ to $SL_h(2)$.
The matrix $\hat{R}_h$ has two eigenvalues ($1$ and $-1$) and a 
spectral decomposition in terms 
of a rank three projector $P_{h+}$ and a rank one projector $P_{h-}$
\begin{equation}\label{1.7.2}
\hat{R}_h  =P_{h+} - P_{h-} \quad, \qquad  P_{h \pm} \hat{R}_h = \pm P_{h \pm} 
\quad, 
\end{equation}
\begin{equation}\label{1.7.3}
P_{h+}= \frac{1}{2}(I+\hat{R}_h) \quad ,\quad P_{h-}= \frac{1}{2} (I-\hat{R}_h)
= \frac{1}{2} \left[ 
\begin{array}{cccc}
0 & h & -h & -h^2 \\
0 & 1 & -1 & -h \\
0 & -1 & 1 & h \\
0 & 0 & 0 & 0 
\end{array} \right]  \quad.
\end{equation}
\noindent
 The deformed determinant $\xi$
in (\ref{1.7}) may be then expressed as
\begin{equation}\label{adet}
(det_h M)P_{h-} := P_{h-} M_{1} M_{2}\quad, \quad
(det_h M^{-1})P_{h-} = M_{2}^{-1} M_{1}^{-1}P_{h-} \quad,
\end{equation}
$$
det_h M^{-1}= (det_h M)^{-1} \quad , \quad (det_h M^{\dagger})P_{h-}^{\dagger}
= M_{2}^{\dagger} M_{1}^{\dagger}P_{h-}^{\dagger} \quad .
$$
\noindent
The following relations  have an obvious  equivalent in the undeformed  
case:
\begin{equation}\label{1.7.4}
\epsilon_hM^t\epsilon_h^{-1}=M^{-1} \,, \, \epsilon_h= \left(
\begin{array}{cc}
      h & 1 \\
      -1 & 0
      \end{array} \right) \,, \,  \epsilon_h^{-1}= \left(
\begin{array}{cc}
      0 & -1 \\
      1 & h
      \end{array} \right) \,, \, 
P_{h- \;ij,kl}=\frac{-1}{2}\epsilon_{h\;ij}\epsilon^{-1}_{h\;kl} \,.
\end{equation}
\noindent

\vspace{0.5\baselineskip}

As in the standard $q$-case \cite{MANIN}, a `quantum $h$-plane' may be 
defined for $GL_h(2)$. The deformed  $h$-plane
associated  with $GL_h(2)$
is the associative algebra generated by 
two elements $(x,y)\equiv X$,  the commutation properties of which
are  given by \cite{EOW,KAR}
$xy=yx+hy^2$ ($R_hX_1X_2=X_2X_1$). These  commutation
relations  are preserved under  the  transformations 
$X'=MX$. This invariance statement,   suitably extended  to apply
to the case of deformed Minkowski 
spaces, provides the essential 
ingredient for a  classification  of the deformations of the Lorentz group
\cite{WOZA} and of the associated Minkowski algebras \cite{JPA96}
(see also \cite{PODWOR}; we shall not consider here deformations governed
by a dimensionful parameter).

 Deformed `groups' related with  
different values of $h \in C$
are equivalent and  their $R_h$ matrices  are related by  a similarity 
transformation.  Thus, from now on, we shall take  
$h \in R$.

\section{$h$-deformed  Lorentz groups}

\indent

 The determination of a complete set of 
deformations of  the Lorentz group (see \cite{WOZA,JPA96}) requires replacing 
(see \cite{POWO,WA-ZPC48,OSWZ-CMP})
the $SL(2,C)$ matrices $A$ in $K'=AKA^{\dagger}$ 
($K=K^{\dagger}=\sigma_{\mu}x^{\mu}$) by the generator 
matrix $M$ of a deformation of $SL(2,C)$,  and 
the characterization of all possible commutation relations
among the generators 
($a,b,c,d$) of $M$ and ($a^*,b^*,c^*,d^*$) of $M^{\dagger}$.

In particular, for  the deformed Lorentz  groups associated with $SL_h(2)$,
the $R$-matrix form of these commutation relations  may be   expressed  by
\begin{equation}\label{2.3}
\begin{array}{ll}
R_h M_1M_2=M_2M_1 R_h \quad, & \quad
M_1^{\dagger}R^{(2)}M_2=M_2R^{(2)}M_1^{\dagger} \quad, \\
M_2^{\dagger}R^{(3)}M_1=M_1R^{(3)}M_2^{\dagger} \quad, & \quad
R^{\dagger}_h M_1^{\dagger}M_2^{\dagger}
=M_2^{\dagger}M_1^{\dagger}R^{\dagger}_h \quad,
\end{array}
\end{equation}
\noindent
where $R^{(3) \,\dagger}=R^{(2)}={\cal P}R^{(3)}{\cal P}$
(`reality' condition for $R^{(3)}$). 
The consistency of these relations is assured if 
$R^{(3)}$ satisfies the `{\it mixed} Yang-Baxter-like'
equation  (see \cite{FTUV94-21,FM} in this respect)
\begin{equation}\label{2.4}
R_{h\, 12}R^{(3)}_{13}R^{(3)}_{23}=R^{(3)}_{23}R^{(3)}_{13}R_{h\,12}
\quad.
\end{equation}
\noindent
This equation, considered as an `FRT' equation, 
indicates that $R^{(3)}$ is a representation of  $GL_h(2)$,  
$(M_{ij})_{\alpha \beta}=
R^{(3)}_{i \alpha , j \beta}$. Thus, $R^{(3)}$ may be seen as a 
matrix in which its 2$\times$2 blocks satisfy  
among themselves the same commutation relations as  the entries of $M$,
\begin{equation}
R^{(3)}= \left[ \begin{array}{cc}
                   A & B \\
                   C & D 
                 \end{array} \right] \; \sim \;
M = \left[ \begin{array}{cc}
                   a & b \\
                   c & d 
                 \end{array} \right] \;.
\end{equation}
\noindent
As a result, the problem of finding all 
possible Lorentz $h$-deformations is
equivalent to finding all possible $R^{(3)}$ matrices with 2$\times$2 block
entries satisfying (\ref{1.6}) such that 
${\cal P}R^{(3)}{\cal P} = R^{(3) \dagger}$ ($ \hat{R}^{(3)}
=\hat{R}^{(3) \dagger}$).

The solutions of these equations are (see  \cite{JPA96,WOZA}) ($h \in R$)
\begin{equation}\label{4.8}
\mbox{{\bf 1.}} \quad  R^{(3)}=\left[
\begin{array}{cccc}
1 & 0 & 0 & 0  \\
0  & 1 & r  & 0   \\
0  & 0 & 1 & 0 \\
0  & 0 & 0 & 1
\end{array} \right] \;, \;   r \in R \,;   \quad
\mbox{{\bf 2.}} \quad   
R^{(3)}=\left[
\begin{array}{cccc}
1 & 0 & -h & 0  \\
-h  & 1 & 0  & h   \\
0  & 0 & 1 & 0 \\
0  & 0 & h & 1
\end{array} \right] \;. 
\end{equation}
\noindent
They characterize the two  $h$-deformed Lorentz groups, which
will be denoted $L_h^{(1)}$ and $L_h^{(2)}$ respectively.
Using  (\ref{4.8})  in (\ref{2.3}), the commutation 
relations among the entries of $M$ and $M^{\dagger}$ read
(see also \cite{WOZA})

\vspace{0.5\baselineskip}
 
\noindent
{\bf 1.}  $L_h^{(1)}$:
\begin{equation}\label{mm1}
\begin{array}{lll}
{} [a,a^*]=rc^*c  \;, & [a,b^*]=rd^*c \;, & [a,d^*]=0 \;, \\
{} [b,b^*]=r(d^*d - aa^*) \;, & [b,d^*]=-rc^*a  \;, & [d,d^*]=-rc^*c \;, \\
{} [c,M^{\dagger}]=0 \;. & \, & \,
\end{array} 
\end{equation}

\noindent
{\bf 2.}  $L_h^{(2)}$:
\begin{equation}\label{mm2}
\begin{array}{ll}
{} [a,a^*]=-h(c^*a + a^*c) \;, & [a,b^*]=h(aa^*-d^*a-b^*c) \;,  \\
{} [a,c^*]= hc^*c \;, & [a,d^*]= h(ac^* + d^*c) \;,  \\
{} [b,b^*]=-h(ab^* + ba^* + b^*d+d^*b) \;, & [b,c^*]= h(c^*d+ac^*) \;, \\
{} [b,d^*]= h( d^*d - ad^*-bc^*) \;, & [c,c^*]=0 \;, \\
{} [c,d^*]= hc^*c \;, & [d,d^*]= -h(cd^*+dc^*) \;.
\end{array}
\end{equation}

\section{$h$-Deformed Minkowski spaces }

\indent

To introduce the deformed Minkowski {\it algebra} ${\cal M}^{(j)}_h$ 
associated with a deformed Lorentz group  $L^{(j)}_h$ ($j$=1,2)
it is natural to extend   $K'=AKA^{\dagger}$ above
by stating that in  the deformed case   the corresponding 
$K$ generates  a comodule algebra for the coaction $\phi$ defined by
\begin{equation}\label{2.5}
\phi : K \longmapsto
K' = M KM^{\dagger} \;, \; K'_{is} = M_{ij} M^{\dagger}_{ls}
K_{jl} \;, \; K=K^{\dagger} \;,\; \Lambda = M \otimes M^* \;,
\end{equation}
\noindent
where, as usual,  it is assumed that the elements of $K$, 
which now do not commute among themselves, commute with those of $M$
and $M^{\dagger}$.  
As in Sec. 1 for $h$-two-vectors (rather, $h$-two-{\it spinors})
we now demand that the commuting 
properties of the entries  of $K$  are preserved by
(\ref{2.5}).  The use of covariance arguments 
to characterize the algebra generated by the elements of  $K$
has been extensively used, and the resulting equations are associated 
with the name of reflection equations \cite{K-SKL,KS} or, in a more
general setting, braided algebras \cite{MAJ-LNM,SMBR2}. In  the present
$SL_h(2)$ case, the 
commutation properties of the entries of the hermitian  matrix $K$ generating 
a deformed Minkowski algebra ${\cal M}_h$ are given by a 
reflection equation of the form 
\begin{equation}\label{2.6}
R_h K_{1} R^{(2)} K_{2} = K_{2} R^{(3)} K_{1} R^{\dagger}_h \quad,
\end{equation}
\noindent
where the $R^{(3)}=R^{(2)\, \dagger}$ matrices are those given in 
(\ref{4.8}).
Indeed, writing equation (\ref{2.6}) for $K'=MKM^{\dagger}$, 
it follows that the invariance of the commutation properties of
$K$  under the associated deformed Lorentz transformation
(\ref{2.5}) is achieved if relations (\ref{2.3}) are satisfied. 

The deformed  Minkowski length and metric, invariant under 
a Lorentz transformation
(\ref{2.5}) of $L^{(j)}_h$, is defined through the quantum determinant of $K$
given by  \cite{JPA96}
\begin{equation}\label{gdet}
(det_{h}K)P_{h-}P_{h-}^{\dagger} 
= - P_{h-} K_1 \hat{R}^{(3)}K_1 P_{h-}^{\dagger} \quad,
\end{equation}
\noindent
where $P_{h-}P_{h-}^{\dagger}$ is a projector since
$(P_{h-}P_{h-}^{\dagger})^2=
\left( \frac{2+h^2}{2} \right)^2 P_{h-}P_{h-}^{\dagger}$.
The above $h$-determinant is invariant, central and, since $\hat{R}^{(3)}$ 
and $K$ are hermitian, real; thus, it defines  the 
{\it deformed Minkowski length} $l_h$
for the $h$-deformed spacetimes ${\cal M}^{(j)}_h$.

Similarly, it is possible to write the $L_h^{(j)}$-invariant scalar product 
of {\it contravariant} (transforming  as the matrix $K$, eq.  (\ref{2.5}))
and {\it covariant} ($Y \mapsto Y'=(M^{\dagger})^{-1}YM^{-1}$)
matrices (Minkowski four-vectors)  as the quantum trace (\cite{FRT1,ZUM})
of a matrix product 
\cite{FTUV94-21,JPA96}.  We define  the 
$h$-deformed trace of a matrix $B$ by 
\begin{equation}\label{gtr}
tr_{h}(B):=tr(D_hB) \quad , \quad  
D_h:=tr_{(2)}({\cal P}(( R_h^{t_1})^{-1})^{t_1})=
\left( \begin{array}{cc} 
        1 & -2h \\
          0   &  1
\end{array} \right) \;, 
\end{equation}
\noindent
where $tr_{(2)}$ means trace in the second space.
This deformed trace is invariant under the quantum group coaction
$B \mapsto MBM^{-1}$ since the expression of $D_h$ above guarantees 
that $D_h^t=M^tD_h^t(M^{-1})^t$ is fulfilled.
To check this explicitly, we start by transposing the first eq. in (\ref{2.3})
in the first space, obtaining 
\begin{equation}\label{-al}
M_{1}^{t_{1}} R_h^{t_{1}} M_{2} = M_{2} R_h^{t_{1}} M_{1}^{t_{1}}\quad.
\end{equation}
\noindent
Inverting this expression and multiplying by $M_{1}^{t_1}$ from left and right,
we get 
\begin{equation}\label{-am1}
(M_{1}^{t_{1}})_{ia} (M_{2}^{-1})_{jb} (R_h^{t_{1}})^{-1}_{ab,kl}= 
(R_h^{t_{1}})_{ij,st}^{-1} (M_{2} ^{-1})_{tl} (M_{1} ^{t_{1}})_{sk}\quad.
\end{equation}
\noindent
Making $l = k$ and summing over $k$ this gives
\begin{equation}\label{-an}
M_{ia}^{t} (M^{-1})_{jb} (R_h^{t_{1}})^{-1}_{ab,kk} = (R_h^{t_{1}})^{-1}_{ij,ss}
\quad , \quad
(R_h^{t_{1}})^{-1}_{ij,ss} = M_{ia}^{t} (R_h^{t_{1}})^{-1}_{ab,kk} 
(M^{-1})^{t}_{bj} \quad .
\end{equation}
This allows us to define $D_h$  as $D_{h\,ij} 
=(R_h^{t_{1}})^{-1}_{ji,ss}$ so that (\ref{gtr}) is obtained.

Let us now find the expression of the metric tensor. Consider
\begin{equation}\label{Reps}
K^{\epsilon}_{ij}:=\hat{R}^{\epsilon}_{h \, ij,kl}K_{kl} \quad, 
\quad  \hat{R}^{\epsilon}_h  \equiv
(1 \otimes (\epsilon_h^{-1})^t) 
\hat{R}^{(3)} (1 \otimes (\epsilon_h^{-1})^{\dagger}) \quad .
\end{equation}
\noindent 
Then, if $K$ is contravariant [(\ref{2.5})],  $K^{\epsilon}$ is covariant
{\it i.e.},
$K^{\epsilon} \mapsto (M^{\dagger})^{-1} K^{\epsilon} M^{-1}$.
This may be checked  by using the property of $\hat{R}^{\epsilon}_h$, 
\begin{equation}\label{ap24}
\hat{R}^{\epsilon}_h (M \otimes (M^{\dagger})^t)
=( (M^{\dagger})^{-1} \otimes (M^{-1})^t ) \hat{R}^{\epsilon}_h 
\; \;  \mbox{or} \;\;  \hat{R}^{\epsilon}_h M_1M^*_2=
(M^{\dagger}_1)^{-1} (M^{-1}_2)^t \hat{R}^{\epsilon}_h \;,
\end{equation}
\noindent
which follows from the preservation of the $h$-symplectic metric $\epsilon_h$.
Now, using the expression of  $\epsilon_h$ in (\ref{1.7.4}),
$(P_{h-})_{ij,kl}=- \frac{1}{2}\epsilon_{h\;ij}\epsilon^{-1}_{h\;kl}$ 
and  $D_h= -\epsilon_h(\epsilon_h^{-1})^t$, it
follows that the {\it $h$-deformed  Minkowski length} $l_h$ and 
{\it $h$-metric} $g_h$ are  given by
\begin{equation}\label{hmetr}
\begin{array}{c}
l_h = det_hK = \frac{1}{2+h^2}tr_hKK^{\epsilon} \equiv
g_{h\, ij,kl} K_{ij}K_{kl} \quad , \quad
g_{h\, ij,kl}= \frac{1}{2+h^2} D_{h\, si}\hat{R}^{\epsilon}_{h \, js,kl}\;. 
\end{array}
\end{equation}
\noindent
The $h$-metric is preserved under   $h$-Lorentz transformations
$\Lambda=M \otimes M^*$,
\begin{equation}\label{ap23}
\begin{array}{c}
\Lambda^t g_h \Lambda=g_h
\quad, \quad  g_h= \frac{1}{ 2+h^2}
(D_h^t \otimes 1) {\cal P} \hat{R}^{\epsilon}_h  \quad .  
\end{array}
\end{equation}
\noindent
This is  checked using   eq. (\ref{ap24}) and that 
$D_h^t=M^t D_h^t (M^{-1})^t$:
\begin{equation}\label{ap26}
\begin{array}{ll}
\Lambda^t g_h \Lambda & = (M \otimes M^*)^t g_h (M \otimes M^*) =
 M^t_1 M^{\dagger}_2 g_h M_1 M^*_2 \\
\, & = \frac{1}{2+h^2}
M^t_1 M^{\dagger}_2 D^t_{h\, 1}{\cal P} \hat{R}^{\epsilon}_h
M_1M^*_2 = \frac{1}{2+h^2}
{\cal P} M^t_2 M^{\dagger}_1 D^t_{h\, 2} \hat{R}^{\epsilon}_h
M_1M^*_2 \\
\, & = \frac{1}{2+h^2}
{\cal P} M^t_2 M^{\dagger}_1 D^t_{h\, 2} (M^{\dagger}_1)^{-1}
(M^{-1}_2)^t \hat{R}^{\epsilon}_h 
= \frac{1}{2+h^2}
{\cal P} M^t_2  D^t_{h\, 2} 
(M^{-1}_2)^t \hat{R}^{\epsilon}_h \\
\, &= \frac{1}{2+h^2}
{\cal P}  D^t_{h\, 2}  \hat{R}^{\epsilon}_h 
= \frac{1}{2+h^2}
D^t_{h\, 1}{\cal P}  \hat{R}^{\epsilon}_h = g_h \quad.
\end{array}
\end{equation}

\subsubsection*{The deformed $h$-Minkowski algebras ${\cal M}_h^{(j)}$} 

\indent

Using (\ref{4.8}) in eqs. (\ref{2.6}), 
(\ref{gdet})
and (\ref{hmetr}), the $h$-Minkowski algebras associated with $SL_h(2)$ as 
well as the deformed Minkowski length and metric read

\vspace{0.5\baselineskip}

\noindent
{\bf 1.} ${\cal M}_h^{(1)}$:  Here, $R^{(3)}$ is the first matrix in
(\ref{4.8}). Then, eq.  (\ref{2.6}) gives ($h$ real)
\begin{equation}\label{5.7}
\begin{array}{ll}
{} [ \alpha, \beta ]= -h \beta^2 -r \beta \delta + h \delta \alpha 
-h \beta \gamma + h^2 \delta \gamma \;,&
\quad [\alpha, \delta]= h( \delta \gamma - \beta \delta) \;,\\
{} [\alpha, \gamma]= h \gamma^2 + r \delta \gamma - h \alpha \delta
+h \beta \gamma - h^2 \beta \delta \;, & \quad 
[\beta, \delta]=  h \delta^2 \;,\\
{} [ \beta, \gamma]= h \delta ( \gamma + \beta) + r \delta^2 \;,& \quad
[ \gamma, \delta]= - h \delta^2 \;;
\end{array}
\end{equation}
\begin{equation}\label{5.7.1}
det_{h} K= \frac{2}{h^2+2}(\alpha \delta - \beta \gamma + h \beta \delta) \quad;
\end{equation}
\small
\begin{equation}\label{reps1}
\hat{R}_h^{\epsilon}=\left[
\begin{array}{cccc}
0 & 0 & 0 & 1  \\
0  & -1 & 0  & h   \\
0  & 0 & -1 & h \\
1  & -h & -h & h^2-r
\end{array} \right] \;, \; 
K^{\epsilon}= \left[ \begin{array}{cc}
                  \delta & - \beta + h \delta \\
                - \gamma + h \delta &  \alpha - h ( \beta + \gamma) + 
                                     (h^2 - r) \delta 
                      \end{array} \right] \;;
\end{equation}
\normalsize
\begin{equation}\label{hmetr1}
g_h=  \frac{1}{2+h^2} \left[
\begin{array}{cccc}
0 & 0 & 0 & 1  \\
0  & 0  & -1  & h  \\
0  & -1 & 0 & -h \\
1  & -h & h & -r-h^2
\end{array} \right] \quad.
\end{equation}

\vspace{0.5\baselineskip}

\noindent
{\bf 2.} ${\cal M}_h^{(2)}$: In this case, $R^{(3)}$ is the second matrix in
(\ref{4.8}). Then,
\begin{equation}\label{5.10}
\begin{array}{ll}
{} [ \alpha, \beta ]= 2h \alpha \delta + h^2 \beta \delta  \quad , &
\quad [\alpha, \delta]= 2h (\delta \gamma -  \beta \delta ) \quad,\\
{} [\alpha, \gamma]= -h^2 \delta \gamma - 2h \delta \alpha \quad, & \quad 
[\beta, \delta]=  2 h \delta^2 \quad , \\
{} [ \beta, \gamma]= 3h^2 \delta^2  \quad ,& \quad
[ \gamma, \delta]= - 2h \delta^2 \quad ;
\end{array}
\end{equation}
\begin{equation}\label{5.11}
det_{h} K=\frac{2}{h^2+2}(\alpha \delta - \beta \gamma + 2h \beta \delta) \quad;
\end{equation}

\begin{equation}\label{reps2}
\hat{R}^{\epsilon}_h=\left[
\begin{array}{cccc}
0 & 0 & 0 & 1  \\
0  & -1 & 0  & 2h   \\
0  & 0 & -1 & 2h \\
1  & 0 & 0 & h^2
\end{array} \right] \quad, \quad 
K^{\epsilon}= \left[ \begin{array}{cc}
                 \delta & -\beta + 2h \delta   \\
                 - \gamma + 2h \delta & \alpha + h^2 \delta 
                      \end{array} \right] \quad;
\end{equation}
\begin{equation}\label{hmetr2}
g_h=\frac{1}{2+h^2} \left[
\begin{array}{cccc}
0 & 0 & 0 & 1  \\
0   & 0 & -1  & 2h   \\
0  & -1 & 0 & 0 \\
1  & 0 & 2h & -3h^2
\end{array} \right] \quad.
\end{equation}

\vspace{0.5\baselineskip}

  We might define a `time' generator for these $h$-deformed spacetimes
as proportional to  $tr_hK$ (=$2x^0$ in the undeformed case). However, 
the resulting  algebra element has undesirable properties:  
$tr_hK$ is neither real and nor central. The time generator is central
only for the $q$-deformed Minkowski space of 
\cite{WA-ZPC48,OSWZ-CMP}
(${\cal M}_q^{(1)}$ in the notation of \cite{FTUV94-21,JPA96}).

\section{Differential calculus on ${\cal M}_h^{(j)}$}

\indent

To describe the differential calculus on $h$-Minkowski spaces,
we need expressing the different commutation relations among the fundamental 
objects: deformed coordinates, derivatives and  one-forms.
Following the approach of \cite{AKR,FTUV94-21,KARPACZ} 
to the differential calculus on 
Minkowski algebras associated with the standard deformation $SL_q(2)$, we 
introduce the reflection equations expressing the  commutation
relations defining the algebras of $h$-derivatives and $h$-one-forms
($h$-differential calculi have been considered in \cite{KAR} and 
in \cite{PODLES}  for quantum $N$-dimensional homogeneous spaces
\footnote{We are grateful to the referee for drawing our attention
to \cite{PODLES}, where similar results were independently obtained.}).
The triangularity property of $R_h$ 
provides  these algebras with some advantages  with respect to the 
$q$-deformed ones; namely, the invariance requirement leads in both
cases to only one  reflection equation. This is due to the fact that, 
in general, the equivalent `FRT' equations 
\begin{equation}\label{1.1a}
R_{12}M_1M_2=M_2M_1R_{12} \quad, \quad 
R_{21}^{-1}M_1M_2=M_2M_1R_{21}^{-1} \quad, 
\end{equation}	
\noindent
allow us to take in $R^{(1)} K_{1} R^{(2)} K_{2} 
= K_{2} R^{(3)} K_{1} R^{(4)}$ (cf. (\ref{2.6})) 
$R^{(1)}=R_{12}$ or $R_{21}^{-1}$,
$R^{(4)}=R_{12}^{\dagger}$ or $(R_{21}^{-1})^{\dagger}$
(see \cite{JPA96}). When the triangularity condition holds, however,
$R_{12}=R_{21}^{-1}$ and there is only one possibility. This argument also 
applies to algebras other than the coordinates algebra. Moreover, we shall show 
that  it is possible to introduce
`coordinates' and `derivatives' which are respectively and {\it simultaneously} 
hermitian and antihermitian. 

\subsubsection*{The algebras of $h$-deformed derivatives ${\cal D}_h^{(j)}$}

\indent

As in \cite{FTUV94-21}, we introduce the derivatives by means of an object 
$Y$ transforming  covariantly {\it i.e.},

\begin{equation}\label{dif1}
Y \longmapsto Y'=(M^{\dagger})^{-1} Y M^{-1} \quad , \quad
Y= \left[ \begin{array}{cc}
                  \partial_\alpha & \partial_\gamma  \\
                 \partial_\beta & \partial_\delta 
                      \end{array} \right] \quad;
\end{equation}
\noindent
The  commutation properties of the derivatives are described by

\begin{equation}\label{dif2}
R_h^{\dagger}Y_2R^{(2) \,-1}Y_1=Y_1R^{(3) \,-1}Y_2R_h \quad,
\end{equation}
\noindent
where $R^{(3)}={\cal P}R^{(2)}{\cal P}$ is given in (\ref{4.8}), 
and are preserved under the $h$-Lorentz coaction.  
Since the covariance requirement is the main
ingredient in our  approach let us check explicitly that
(\ref{dif2})  is invariant under 
(\ref{dif1}). Multiplying (\ref{dif2}) by 
$(M_1^{\dagger})^{-1}(M_2^{\dagger})^{-1}$  from the left
and by $M_2^{-1}M_1^{-1}$ from the right 
we get,  using the first and the
last equations in (\ref{2.3}),
\small
\begin{equation}
R_h^{\dagger} (M_2^{\dagger})^{-1} Y_2 (M_1^{\dagger})^{-1} R^{(2) \,-1}
M_2^{-1} Y_1 M_1^{-1} = (M_1^{\dagger})^{-1} Y_1 (M_2^{\dagger})^{-1} 
R^{(3) \,-1} M_1^{-1} Y_2 M_2^{-1} R_h \;.
\end{equation}
\normalsize
\noindent
Finally, using the second and third third eqs. in (\ref{2.3}), we obtain 
\small
\begin{equation}
R_h^{\dagger} \, (M_2^{\dagger})^{-1} Y_2 M_2^{-1}\, R^{(2) \,-1}\,
(M_1^{\dagger})^{-1} Y_1 M_1^{-1} = (M_1^{\dagger})^{-1} Y_1  M_1^{-1} \,
R^{(3) \,-1}\,(M_2^{\dagger})^{-1}  Y_2 M_2^{-1} R_h \;,
\end{equation}
\normalsize
\noindent
which is eq. (\ref{dif2}) for $Y'=(M^{\dagger})^{-1}YM^{-1}$.

The $h$-deformed d'Alembertian may be  introduced by using the $h$-trace
\begin{equation}\label{dif3}
\Box_h \equiv  \frac{1}{2+h^2} tr_h(Y^{\epsilon}Y) \quad, \quad 
Y^{\epsilon}= (\hat{R}_h^{\epsilon})^{-1} Y \quad.
\end{equation}
\noindent
As $l_h$,  $\Box_h$  is Lorentz invariant, real and central 
in the algebra ${\cal D}_h^{(j)}$ of  derivatives.

\vspace{0.5\baselineskip}

Using (\ref{4.8}) in eqs. (\ref{dif2}), 
the commutation relations for ${\cal D}_h^{(j)}$ read

\vspace{0.5\baselineskip}

\noindent
{\bf 1.} ${\cal D}_h^{(1)}$:
\small
\begin{equation}\label{dh1}
\begin{array}{ll}
{} [ \partial_{\alpha} , \partial_{\beta}] = -h \partial_{\alpha}^2 
 \;,&
[\partial_{\beta}, \partial_{\delta}]=  h ( \partial_{\beta}^2
 + \partial_{\beta} \partial_{\gamma} - \partial_{\alpha}
\partial_{\delta}) -r \partial_{\beta} \partial_{\alpha} - h^2 \partial_{\alpha}
\partial_{\gamma} \;,\\
{} [\partial_{\alpha}, \partial_{\gamma}]= h \partial_{\alpha}^2
 \;, & 
[ \partial_{\gamma}, \partial_{\delta}]= -h (\partial_{\gamma}^2
 +  \partial_{\beta} \partial_{\gamma} -  \partial_{\delta}
\partial_{\alpha}) + r \partial_{\alpha} \partial_{\gamma} + h^2 \partial_{\beta}
\partial_{\alpha} \;, \\
{} [\partial_{\alpha}, \partial_{\delta}]= h(
\partial_{\alpha} \partial_{\beta} - \partial_{\gamma} \partial_{\alpha})
 \;, &  
 [ \partial_{\beta}, \partial_{\gamma}]= h \partial_{\alpha} (
\partial_{\beta} + \partial_{\gamma}) -r \partial_{\alpha}^2 \;; 
\end{array}
\end{equation}
\normalsize

\noindent
{\bf 2.} ${\cal D}_h^{(2)}$:
\begin{equation}\label{dh2}
\begin{array}{ll}
{} [ \partial_{\alpha} , \partial_{\beta}] = - 2h \partial_{\alpha}^2 
 \quad ,& \quad 
[\partial_{\alpha}, \partial_{\delta}]= 2h(
\partial_{\beta} \partial_{\alpha} - \partial_{\alpha} \partial_{\gamma})
 \;, \\
{} [\partial_{\alpha}, \partial_{\gamma}]= 2 h \partial_{\alpha}^2
 \quad , & \quad
[\partial_{\beta}, \partial_{\delta}]=  h^2 \partial_{\beta} 
\partial_{\alpha} -2h   \partial_{\delta} \partial_{\alpha}  \;, \\
{} [ \partial_{\beta}, \partial_{\gamma}]= 5h^2 \partial_{\alpha}^2 \quad ,
& \quad 
[ \partial_{\gamma}, \partial_{\delta}]= 
2h \partial_{\alpha} \partial_{\delta} - h^2 \partial_{\alpha} 
\partial_{\gamma}  \;.  
\end{array}
\end{equation}

\subsubsection*{Commutation relations for coordinates and derivatives}

\indent

The commutation relations among the entries of $K$ and $Y$ may be expressed 
by an inhomogeneous reflection equation (see \cite{AKR,FTUV94-21})
\begin{equation}\label{dif4}
Y_2R_hK_1R^{(2)}=R^{(3)}K_1R^{\dagger}_hY_2 + \eta R^{(3)}{\cal P} \quad,
\end{equation}
\noindent
which extends to the $h$-deformed case the undeformed  relation
$\partial_{\mu}x^{\nu}= \delta^{\nu}_{\mu} + x^{\nu}\partial_{\mu}$.
This equation is  consistent with the commutation relations defining the 
algebras  ${\cal M}_h^{(j)}$, ${\cal D}_h^{(j)}$,  and  
is  invariant under $h$-Lorentz transformations   (as already mentioned, 
there is only one $L_h$-invariant reflection equation
 due to the triangularity property
of $R_h$). The invariance  is  seen by multiplying eq. (\ref{dif4})
by $(M_2^{\dagger})^{-1}M_1$  from the left and by $M_1^{\dagger}M_2^{-1}$
from the right and using the commutation relations in (\ref{2.3}).

\vspace{0.5\baselineskip}

It is a common feature of all $q$-deformed Minkowski spaces that the 
covariance transformation properties for `coordinates' and
`derivatives' are consistent with their hermiticity. The mixed commutation 
relations (as expressed by an  inhomogeneous reflection 
equation), however, do not allow in general for 
{\it  simultaneously} hermitian  coordinates and derivatives, 
a feature of non-commutative geometry 
already noted in  \cite{OSWZ-CMP,OZ}.
Let us then look at  the hermiticity properties of $K$ and $Y$ for our
$h$-deformed Minkowski spaces. Clearly, eqs. (\ref{2.6}) and (\ref{dif2}) 
allow us to take both $K$ and $Y$ hermitian.
Keeping the physically reasonable assumption that $K$ is hermitian, eq.
(\ref{dif4}) gives
\begin{equation}\label{dif5}
R^{(2) \dagger } K_1  R^{\dagger}_hY_2^{\dagger}=
Y_2^{\dagger}R_hK_1R^{(3) \dagger} + \eta^* {\cal P} R^{(3) \dagger} \quad.
\end{equation}
\noindent
Since $R^{(3) \dagger}=R^{(2)}={\cal P}R^{(3)} {\cal P}$, we get that
$(-Y^{\dagger})$ satisfies (for $\eta^*=\eta$) the same 
commutation relations than $Y$. Thus, eqs. (\ref{2.6}), (\ref{dif1})
and (\ref{dif4}) are compatible with the  hermiticity of $K$ and the 
antihermiticity  of $Y$. This linear conjugation structure, absent in 
the $q$-deformation, may facilitate the formulation of invariant 
field equations on  ${\cal M}_h^{(j)}$.

\vspace{0.5\baselineskip}

Using the $R^{(3)}$ matrices in (\ref{4.8}) in eq. (\ref{dif4}) 
and setting $\eta =1$, the mixed commutation relations are found to be
(we give only a few cases)

\vspace{0.5\baselineskip}

\noindent
{\bf 1.} ${\cal M}_h^{(1)} \times {\cal D}_h^{(1)}$:
\begin{equation}\label{KY1}
\begin{array}{ll}
{} [\partial_{\alpha}, \alpha]=
1+ h \beta \partial_{\alpha} - h \gamma \partial_{\alpha} - h^2 \delta 
\partial_{\alpha} \quad , & \quad 
 [\partial_{\alpha}, \beta ] = -h \delta \partial_{\alpha} \quad, \\
{} [\partial_{\alpha}, \gamma ] = h \delta \partial_{\alpha} \quad , & \quad
[\partial_{\alpha}, \delta ] = 0 \quad , 
\end{array}   
\end{equation}
$$
\begin{array}{l}
{} [\partial_{\beta}, \alpha]=
(r+h^2) \gamma \partial_{\alpha} + h (r- h^2) \delta \partial_{\alpha} 
- h( \alpha \partial_{\alpha} + \beta \partial_{\beta} 
+ \gamma \partial_{\beta}) +h^2 (\delta \partial_{\beta}  
+  \beta \partial_{\alpha}) \;, \\ 
{} [\partial_{\beta}, \beta ] = 1+ (r- h^2) \delta \partial_{\alpha} 
+ h( \beta \partial_{\alpha} -  \delta \partial_{\beta}) \;, \\
{} [\partial_{\beta}, \gamma ] = 
- h( \gamma \partial_{\alpha} + \delta \partial_{\beta})
+h^2 \delta \partial_{\alpha} \;, \\
{} [\partial_{\beta}, \delta ] = h \delta \partial_{\alpha} \quad . 
\end{array}
$$

\noindent
{\bf 2.} ${\cal M}_h^{(2)} \times {\cal D}_h^{(2)}$:
\begin{equation}\label{KY2}
\begin{array}{ll}
{} [\partial_{\alpha}, \alpha]=
1+ 2h (\beta \partial_{\alpha} -  \gamma \partial_{\alpha}) -4 h^2 \delta 
\partial_{\alpha} \;, & \quad 
[\partial_{\alpha}, \beta ] = -2h \delta \partial_{\alpha} \;, \\
{} [\partial_{\alpha}, \gamma ] =2 h \delta \partial_{\alpha} \;, &
\quad  [\partial_{\alpha}, \delta ] = 0 \;, \\
\; & \;  \\
{} [\partial_{\beta}, \alpha]= - 2h \alpha \partial_{\alpha} 
-h^2 \gamma \partial_{\alpha} - 2h^3 \delta \partial_{\alpha} 
 \;, & \quad 
[\partial_{\beta}, \delta ] = 2h \delta \partial_{\alpha} \;, \\
{} [\partial_{\beta}, \beta ] = 1 - h^2 \delta \partial_{\alpha} \;, &
\quad  [\partial_{\beta}, \gamma ] = 4h^2 \delta \partial_{\alpha} \;
\end{array}
\end{equation}

\subsubsection*{The algebras of  $h$-deformed one-forms $\Lambda_h^{(j)}$}

\indent

To determine now the commutation relations for the $h$-de Rham complex we
now introduce the exterior derivative $d$ following \cite{OSWZ-CMP}
(see also \cite{FTUV94-21,AKR}). The algebra of the $h$-forms is generated by
the entries of a matrix $dK$. Clearly, $d$ commutes with the 
Lorentz coaction, so that
\begin{equation}\label{dif8}
dK'=M\, dK \, M^{\dagger} \quad.
\end{equation}
\noindent
Applying $d$ to  eq. (\ref{2.6}) we obtain
\begin{equation}\label{dif9}
R_h dK_{1} R^{(2)} K_{2} + R_h K_{1} R^{(2)} dK_{2} = 
dK_{2} R^{(3)} K_{1} R^{\dagger}_h  + K_{2} R^{(3)} dK_{1} R^{\dagger}_h \quad.
\end{equation}
\noindent
Its  only solution is given by
\begin{equation}\label{dif10}
R_h dK_{1} R^{(2)} K_{2} =  K_{2} R^{(3)} dK_{1} R^{\dagger}_h \quad
\end{equation}
\noindent
(which implies  $R_h K_{1} R^{(2)} dK_{2} = 
dK_{2} R^{(3)} K_{1} R^{\dagger}_h$). From eq. (\ref{dif10}), it follows that 
\begin{equation}\label{dif11}
R_h dK_{1} R^{(2)} dK_{2} = - dK_{2} R^{(3)} dK_{1} R^{\dagger}_h \quad.
\end{equation}
\noindent
Again, it is easy to check that these relations are invariant under 
hermitian conjugation. Notice that the reflection equations 
(\ref{2.6}), (\ref{dif10})
and (\ref{dif11})  have the same $R$-matrix structure. 
In the $h$-deformed case,
the reflection equation giving the commutation 
relations among the  generators of two differential algebras is determined 
only by the transformation (covariant or contravariant) law of these 
generators. Thus, there are only three types of reflection
equations, those of (\ref{2.6}), (\ref{dif2}) and
(\ref{dif4}),  as a consequence  of the triangularity  of $SL_h(2)$.
In contrast,  in the $q$-deformation (based on $SL_q(2)$),  
the number of  reflection equation types is larger.

The exterior derivative is given for the two $h$-Minkowski algebras by 
\begin{equation}\label{dif12}
d= tr_h(dK\,Y)= d\alpha \partial_\alpha  + d\beta \partial_\beta  +
d\gamma \partial_\gamma  + d\delta \partial_\delta -2h(
d\gamma \partial_\alpha  + d\delta \partial_\beta ) \;.  
\end{equation}
\noindent

\vspace{0.5\baselineskip}

To conclude, let us mention that the additive braided group 
\cite{MAJ-LNM,SMBR2}
structure of the above  algebras  may be easily found. It suffices to impose 
{\it e.g.} that 
eq. (\ref{2.6})  is also satisfied by the sum $K+K'$ of two {\it copies}
$K$ and
$K'$ of the given $h$-Minkowski algebra.  This leads to an  
equation of the same type as  (\ref{2.6}) (see the comment above)
\begin{equation}\label{dif13}
R_h K'_{1} R^{(2)} K_{2} = K_{2} R^{(3)} K'_{1} R^{\dagger}_h \quad.
\end{equation}
\noindent
Eq. (\ref{dif13}) is clearly preserved by (\ref{2.5}).
The above discussion could be extended easily to obtain a braided 
differential calculus. The unified braided structure of all $h$-
(and $q$-, we note in passing) deformed Minkowski spaces 
was given in \cite{JPA96} (see 
\cite{MeyerM} for the particular  $q$-deformed case of \cite{OSWZ-CMP}). 
Given the  generality of our presentation,
it is a trivial exercise to introduce a {\it unified} additive\footnote{
Although the formalism allows us to introduce multiplicative braiding, we
wish to point out that the multiplication is not consistent with covariance 
for  Minkowski spaces.}  braided differential 
calculus valid for all the  $h$- (or $q$-) Minkowski algebras. For the case
of the $q$-Minkowski space in  \cite{OSWZ-CMP} we refer to
\cite{VLADI}.

\end{document}